\documentclass[12pt]{amsart}

\usepackage{amsmath,amsthm,amssymb}
\usepackage[T1]{fontenc}
\usepackage{latexsym}
\usepackage{subfigure, epsfig, epsf}
\usepackage{color,fancyhdr,lastpage,graphicx}
\usepackage{xspace}
\usepackage{setspace}
\usepackage{bbm}

\newcommand{\NN}{\mathbb{N}}  \newcommand{\ZZ}{\mathbb{Z}}    
\newcommand{\RR}{\mathbb{R}}    \newcommand{\LL}{\mathbb{L}}

   \newcommand{\mcD}{\mathcal{D}} 
\newcommand{\mcC}{\mathcal{C}}    \newcommand{\mcS}{\mathcal{S}}
\newcommand{\mcH}{\mathcal{H}} \newcommand{\mcA}{\mathcal{A}}

\newcommand{\EE}{\mathbb{E}}

\newcommand{\rv}{random variable }
    \newcommand{\rvs}{random variables }

\newcommand{\as}{almost-surely }

\newcommand{\st}{such that }

\newcommand{\ep}{\epsilon}

\renewcommand{\leq}{\leqslant}
\renewcommand{\geq}{\geqslant}

\newcommand{\al}{\alpha}

\newcommand{\ssk}{\smallskip}
\newcommand{\wrt}{with respect to }
\newcommand{\noi}{\noindent }

\newtheorem{thm}{\hspace{-0.15cm} {\sc Theorem} }

\newtheorem{prop}[thm]{\hspace{-0.18cm} {\sc Proposition}}

\numberwithin{equation}{section} 

\newenvironment{Dem}{%
    \begin{list}{\hspace{0.6cm}{\sc Proof --}}{%
        \setlength{\topsep}{0pt}%
        \setlength{\leftmargin}{0pt}%
        \setlength{\rightmargin}{0pt}%
        \setlength{\listparindent}{0pt}%
        \setlength{\itemindent}{0pt}%
        \setlength{\parsep}{0pt}%
        \addtolength{\leftmargin}{20pt}%
        \addtolength{\rightmargin}{0pt}%
    } \item }{\hfill{\space $\rhd$}\end{list}\smallskip}

\title{{\sc KPZ in a multidimensional random geometry of multiplicative cascades}}
\date{\today}
\author[{\sc I.F. Bailleul}]{{\sc I.F. Bailleul}}
\address{Statistical Laboratory, Center for Mathematical Sciences, Wilberforce Road, Cambridge, CB3 0WB, UK}
\email{i.bailleul@statslab.cam.ac.uk}
\urladdr{http://www.statslab.cam.ac.uk/~ismael/}

\begin{document}

\begin{abstract}
We show in this note how the one-dimensional KZP formula obtained by Benjamini and Schramm in \cite{BenjaminiSchramm} can be extended to a multidimensional setting.
\end{abstract}

\maketitle

\section{Hausdorff dimension in a nested measure space}
\label{HausdorffDimRevisited}

\subsection{Dimension}

Let $(S,\mcS,\mu)$ be a measure space and suppose given a nested family of countable $\sigma$-algebras $\mcS_n=\sigma\bigl(A^i_n\,;\,i\geq 1\bigr)$, with $A^i_n\in\mcS$ disjoint up to $\mu$ null sets, and $\mu(A^i_n)>0$ for each $i\geq 1$ and $n\geq 1$. Suppose further that $\ep_n:=\sup\,\mu(A^i_n)$ decreases to $0$ as $n$ goes to infinity. Given $s\geq 0$ and $\delta>0$, set for any measurable $E\in\mcS$
$$
\mcH^s_\delta(E) = \inf\,\sum\mu\bigl(A^{i_\al}_{n_\al}\bigr)^s,
$$
where the infimum is over the set of coverings $E\subset\bigcup_{\al\in\mcA}A^{i_\al}_{n_\al}$ of $E$, indexed by a subset $\mcA$ of $\NN^*\times\NN^*$, and \st $\ep_{n_\al}\leq \delta$ for all $\al\in\mcA$. The quantity $\mcH^s_\delta(K)$ increases as $\delta$ decreases to $0$. Set 
$$
\mcH^s(E) = \underset{\delta\downarrow 0}{\lim} \mcH^s_\delta(E).
$$
Like in the usual definition of the Hausdorff dimension of a set, it is easy to see that if 
\begin{itemize}
   \item $\mcH^{s_0}(E)<\infty$ then $\mcH^t(K)=0$ for any $s_0<t$,
   \item $\mcH^{s_0}(E)=\infty$ then $\mcH^s(K)=0$ for any $s<s_0$,
\end{itemize}
so it makes sense to define the dimension $\zeta_\mu(E)$ of $E$ as $\sup\,\{s\geq 0\,;\,\mcH^s(E)=\infty\} = \inf\{t\geq 0\,;\,\mcH^t(E)=0\}$. \textit{As $\mcH^1$ coincides with $\mu$, it follows that $\zeta_\mu(E)\leq 1$, for any $E\in\mcS$.} So only sets with null $\mu$-measure have a dimension smaller than $1$.

\medskip

\noi \textbf{Open question.} Let us work in the space $S = \mcC\bigl([0,1],\RR\bigr)$, with its Borel $\sigma$-algebra and Wiener measure. Define $A_n^{(j,k)}$ as $\Big\{\omega\in\mcC\bigl([0,1],\RR\bigr)\,;\,\omega\bigl((j+1)2^{-n}\bigr)-\omega\bigl(j2^{-n}\bigr)\in \bigl[k2^{-n},(k+1)2^{-n}\bigr)\Big\}$, for $0\leq j\leq 2^n-1$ and $k\in\ZZ$, and set $\mcS_n = \sigma\bigl(A_n^{(j,k)}\,;\,0\leq j\leq 2^n-1,\;k\in\ZZ\bigr)$. Let us call \textit{\textbf{Wiener-Hausdorff dimension}} the above dimension of a measurable subset of $\mcC\bigl([0,1],\RR\bigr)$. \textit{Compute the Wiener-Hausdorff dimension of the set of $\al$-H\"older continuous paths, for $\al\geq \frac{1}{2}$.}

\subsection{Frostman lemma}

If $(S,\mcS)$ is $\RR^d$ with its Borel $\sigma$-algebra, and $\mcS_n$ is the $\sigma$-algebra generated by the dyadic cubes of side $2^{-n}$, then the above definition of dimension coincides with the usual Hausdorff dimension, up to a multiplicative constant $\frac{1}{d}$; see section 2.4, Chap. 2, in \cite{Falconer}. We adopt the above definition of dimension for the sequel. Like its classical counterpart, the above set function $\mcH^s(\cdot)$ can be shown to be an $\bigl(\RR_+\cup\{\infty\}\bigr)$-valued measure on $\bigl(\RR^d,{\sf Bor}(\RR^d)\bigr)$. The Euclidean background will not appear anymore except under the form of the nested family $(\mcS_n)_{n\geq 0}$.

\medskip

Given two points $x,y\in\RR^d$, define the \textit{ball} $B(x,y)$ as the smallest dyadic cube containg $x$ and $y$, and define their ``distance'' as $\mu\bigl(B(x,y)\bigr)$. Define accordingly the ball $B_r(x) = \{y\in\RR^d\,;\,\mu\bigl(B(x,y)\bigr)\leq r\}$. Working exactly as in theorem 4.10 and proposition 4.11 in \cite{Falconer}, one can prove the following proposition. 

\begin{prop}
For any Borel set $E$ with $0<\mcH^s(E)<\infty$, there exists a constant $c$ and a compact set $K\subset E$ with $\mcH^s(K)>0$ \st 
$$
\mcH^s\bigl(K\cap B_r(x)\bigr) \leq cr^s
$$
for all $x\in\RR^d$ and $r>0$.
\end{prop}
It follows classically that the following version of Frostman lemma holds in our setting. Given any non-negative measure $\nu$ on $\bigl(\RR^d,{\sf Bor}(\RR^d)\bigr)$, define its $s$-energy as
$$
I_s(\nu) = \int\hspace{-0.2cm}\int\frac{\nu(dx)\nu(dy)}{\mu\bigl(B(x,y)\bigr)^s}.
$$

\begin{thm}
\label{Frostman}
If $E$ is a Borel set with $0<\mcH^s(E)$, then there exists a non-negative measure $\nu$ with support in (a compact subset of) $E$ \st $I_t(\nu)<\infty$, for all $t<s$. This is in particular the case if $s<\zeta_\mu(E)$.
\end{thm}

\medskip

\noi \textbf{Remark.} The work \cite{RhodesVargas} contains in section 5.1 a similar, though different, notion of dimension in a metric measure space.

\section{A dimension-free KPZ formula}

Let $\mcD_n = \bigcup_{k=1}^{2^{dn}}A^n_k$ be the dyadic ``partition'' of the unit cube of $\RR^d$ by closed dyadic cubes of side length $2^{-n}$. Given $m<n$, each $A^n_k$ is a subset of a unique $A^m_{k(m)}$. Let $W$ be a positive real-valued \rv with $\EE[W]=1$, and let $\bigl\{(W^n_i)_{i=1}^{2^{nd}}\,;\,n\geq 1\bigr\}$ be an iid sequence of \rvs with common law the law of $W$. Define the measure $\mu_n$ by its density $w_n(x)$ \wrt Lebesgue measure. It is constant, equal to $\prod_{m=0}^nW^m_{k(m)}$, on each $A^n_k$. We adopt as in \cite{BenjaminiSchramm} the notation $\ell$ for $\mu\bigl([0,1]^d\bigr)$. It has expectation no greater than $1$.

\begin{prop}
Almost-surely, the measures $\mu_n$ converge weakly to some random measure $\mu$, which does not charge any dyadic hyperplane. It is \as non-null if $\,\EE[W\log W]<d$.
\end{prop} 

\begin{Dem}
The proof works exactly as in the $1$-dimensional proof, with $2^d$ independent copies of $\ell$ rather than only two.
\end{Dem}

The next result generalizes Benjamini and Schramm's result \cite{BenjaminiSchramm} obtained in a one-dimensional setting.

\begin{thm}
Let $E$ be any Borel set of $[0,1]^d$. Denote by $\zeta_0$ its dimension as defined above using Lebesgue measure, and let $\eta$ be its dimension using the random measure $\mu$. Suppose that $\EE[W\log W]<d$, and $\EE[W^{-s}]<\infty$, for all $s\in[0,1)$. Then $\zeta$ is \as a constant and satisfies the identity
$$
2^{\zeta_0}= \frac{2^\zeta}{\EE[W^\zeta]}
$$
\end{thm}

The above conditions are satisfied by an exponential of Gaussian with a small enough variance.

\ssk

\begin{Dem}
The proof mimicks word by word the proof of \cite{BenjaminiSchramm}. Write $|A|$ for the Lebesgue measure of a Borel set $A$. Set, for $s\in [0,1]$, $\phi(s)=s-\ln_2\EE[W^s]$. Note that since the notion of dimension introduced in section \ref{HausdorffDimRevisited} is no greater than $1$ the function $\phi$ is an increasing homeomorphism from $[0,1]$ to itself.

\ssk

\noi \textbf{a)} Lemma 3.3 becomes here: $\EE\bigl[\mu\bigl(B(x,y)\bigr)^s\bigr]\leq \big|B(x,y)\big|^{\phi(s)}$, for all $x,y\in [0,1]^d$.

{\small Note that the balls $B(x,y)$ are always dyadic balls; suppose the given ball belongs to $\mcD_n$, so $\big|B(x,y)\big| = 2^{-nd}$. Then, we have by the independence in the construction of $\mu$
$$
\EE\bigl[\mu\bigl(B(x,y)\bigr)^s\bigr] = 2^{-nd}\,\EE[W^s]^{nd}\,\EE[\ell^s] \leq \{2^{-nd}\}^{\phi(s)} = \big|B(x,y)\big|^{\phi(s)},
$$
as $0\leq s\leq 1$, so $\EE[\ell^s]\leq \EE[\ell]^s = 1$.} It follows directly that we have \as $\phi(\zeta)\leq \zeta_0$.

\ssk

\noi \textbf{b)} The proof that $\phi(\zeta)\geq\zeta_0$, theorem 3.5,  works identically, replacing the usual energy of a measure by its above modification, and using the version of Frostman lemma provided in theorem \ref{Frostman}. A straightforward adaptation of the proof that $\EE[\ell^{-s}]<\infty$ if $\EE[W^{-s}]<\infty$, given in \cite{BenjaminiSchramm}, gives the same result in our setting. Note also that a different choice of H\"older coefficient is needed to prove that the sequence $\nu_n\bigl([0,1]\bigr)$ is uniformly bounded in some $\LL^p$. 
\end{Dem}

Note that the above theorem does not come as a surprise and should actually hold on much more general state spaces than $[0,1]^d$. It should be interesting in particular to investigate what happens on random trees like Galton-Watson trees, and tree-like objects like random fractals.

\end{document}